\documentclass{article}
\usepackage{a4wide}

\usepackage[czech]{babel}

\usepackage{graphicx}
\usepackage[T1]{fontenc}
\usepackage[utf8]{inputenc}

\usepackage{amsmath, amssymb, amsthm, amsfonts, dsfont, bbm}
\usepackage{enumerate}
\usepackage{caption}
\usepackage{stmaryrd}
\usepackage{hyperref,url}

\newcommand{\atm}{\ensuremath{\mathbf{AST}}}
\newcommand{\zf}{\ensuremath{\mathbf{ZF}}}
\newcommand{\zfc}{\ensuremath{\mathbf{ZFC}}}
\newcommand{\zffin}{\ensuremath{\mathbf{ZF}_{\mathrm{fin}}}}

\newcommand{\ist}{\ensuremath{\mathbf{IST}}}

\newcommand{\FN}{\ensuremath{\mathrm{FN}}}
\newcommand{\N}{\ensuremath{\mathrm{N}}}

\newcommand{\Fin}{\ensuremath{\mathrm{Fin}}}

\begin{document}

\title{Vopěnka's Alternative Set Theory \\
in the Mathematical Canon of the 20th Century: \\
Author's Translation from Czech\footnote{This is the author's 
translation of her paper published originally in Czech: Zuzana Haniková, Vopěnkova Alternativní teorie množin v matematickém kánonu 20.~století,
\emph{Filosofický časopis} 70(3), 485--504, 2022. {\tt https://doi.org/10.46854/fc.2022.3r.485} } }
\author{Zuzana Haniková\\
Institute of Computer Science of the Czech Academy of Sciences\\
hanikova@cs.cas.cz}
\date{November 20, 2022}
\maketitle


\renewcommand{\abstractname}{Abstract}

\begin{abstract}
Vopěnka’s Alternative Set Theory can be viewed both as an evolution and as a revolution: 
it is based on his previous experience with nonstandard universes, inspired
by Skolem’s construction of a nonstandard model of arithmetic, and its inception
has been explicitly mentioned as an attempt to axiomatize Robinson’s Nonstandard
Analysis. Vopěnka preferred working in an axiomatic theory to investigating its individual models; 
he also viewed other areas of non\-classical mathematics through this
prism. This article is a contribution to the mapping of the mathematical neighbourhood 
of the Alternative Set Theory, and at the same time, it submits a challenge to
analyze in more detail the genesis and structure of the philosophical links that eventually 
influenced the Alternative Set Theory.
\end{abstract}

\section{Why come back to the Alternative Set Theory?}
\label{sct:procatm}

By the Alternative Set Theory ($\atm$)\footnote{The attribute ``Vopěnka's'' will be omitted
wherever no confusion can arise; likewise, except in titles and quotations, we refrain from capitalizing the name of the theory.} 
we mean the exposition within two monographs by P.~Vopěnka: 
\emph{Mathematics in the Alternative Set Theory} \cite{Vopenka:Teubner}, published in 1979
in an English translation by P.~Hájek, and 
 \emph{Úvod do matematiky v alternatívnej teórii množín} [An Introduction to Mathematics in the Alternative Set Theory] \cite{Vopenka:Alfa},
published in Slovak ten years later in a translation by P.~Zlatoš, 
where the slightly expanded mathematical material is complemented with plenty of motivating,
philosophical and historical content.
Neither of these monographs aspired to provide the full body of then extant mathematical material on the $\atm$.
Such a body would include voluminous series of journal papers, mostly available via the \emph{Czech Digital
Mathematical Library}\footnote{www.dml.cz}, some of which left their mark in Vopěnka's monographs\footnote{According to
\cite{Vopenka:Alfa}, contributions to the development of the $\atm$ were made by
 K.~Čuda, J.~Mlček, A.~Sochor, A.~Vencovská, J.~Chudáček, M.~Resl, K.~Trlifajová, B.~Vojtášková, J.~Sgall, J.~Witzany, 
J.~Guričan, M.~Kalina and P.~Zlatoš.},
as well as the \emph{Proceedings of the 1st Symposium Mathematics in the Alternative Set Theory} \cite{StaraLesna:sympAST}, 
which provided an opportunity to include comprehensive new material on the theory and can thus be taken as a compendium of
results and references thereto.

All of the above needs to be distinguished from the earlier 1972  monograph 
\emph{The Theory of Semisets} \cite{Vopenka-Hajek:Semisets} by Vopěnka and Hájek,  
and the subsequent \emph{Nová infinitní matematika} [New Infinitary Mathematics] \cite{Vopenka:NIM} by Vopěnka, 
published in several volumes in 2015.

As things stand, Vopěnka's lifelong work and contribution have not been systematically assessed and placed within appropriate context in the history and philosophy of mathematics. Moreover, no rigorous scientific biography is available, even if there is a quite informative
text by A.~Sochor  \cite{Sochor:Vopenka} from a 2001 special issue of the \emph{Annals of Pure and Applied Logic} dedicated to Vopěnka, where Sochor acted also as a guest editor. It offers a periodization of Vopěnka's work up to its publication and side by side with the period of work on the $\atm$, it discusses and highlights also the preceding one, remarkable for the important results of the Prague set theory seminar and the publication of \emph{The Theory of Semisets}.
   
Vopěnka continued the research line of his  monographs on the $\atm$ with multiple lectures and journal papers\footnote{
For example \cite{Vopenka:ocojdevatm}, published previously in English as \cite{Vopenka:ASTallabout} in the translation of A.~Vencovská,
and also as \cite{Vopenka:AstPhil1991}.},
and also with comprehensive texts on the history of set theory \cite{Vopenka:vyprNovobarokniMat}, 
where circumstances of the inception of the $\atm$ are considered.
Czech logicians keep reflecting  Vopěnka's works;
recently, there has been Švejdar's paper \cite{Svejdar:VopenkaHajek} or
a monothematic block ``Setkání s Petrem Vopěnkou'' [An Encounter with Petr Vopěnka] \cite{TemFilCasopis2016Vop}  in \emph{Filosofický časopis} 
[The Philosophical Journal], edited by K.~Trlifajová, with contributions reflecting the period
when Vopěnka worked at the Faculty of Mathematics and Physics of Charles University in Prague, particularly while leading his
$\atm$ research group, as well as some newer directions pursued during the time he worked at the Faculty of Arts of the University of
West Bohemia in Pilsen. In the same recent period, several informal texts were published, such as ``Poezie matematiky Petra Vopěnky''
[The Poetry of Petr Vopěnka's Mathematics] \cite{Trlifajova:poezieMatVopenky}, supplemented by excerpts from \emph{Nová infinitní matematika},
or Vopěnka's obituaries by J.~Mlček  \cite{Mlcek:nekrologVopenka} and J.~Chudáček  \cite{Chudacek:nekrologVopenka}.

Some parts of Vopěnka's works on the $\atm$, and some related works, encourage the reader to
perceive the mathematical world in a dichotomous fashion: on the one side, there is Vopěnka's $\atm$,
on the other, what is termed \emph{Cantor's set theory}---thriving yet self-centered and detached from everyday experience.
The notion is quite common in Vopěnka's work, so it seems natural to pause and ask: what is Cantor's set theory?
Sochor's paper \cite{Sochor:AST76} gives a plain definition: it is a \emph{theory of infinity}, proposed by G.~Cantor, 
which makes it possible to study and classify the notion. The definition is intended as comprehensive: namely, it
exceeds the usual understanding of the notion  \emph{classical set theory}, but this is not clear from the term itself
and it is easy to mistake the one for the other.
The mathematical community is portrayed as quite homogeneous in its pragmatic leaning
toward this theory as foundations of mathematics. Vopěnka does discuss plurality of viewpoints in modern foundations of mathematics,
but he does so almost invariably in retrospect, that is, as issues already dealt with. Both of Vopěnka's monographs 
also refer to a crisis in contemporary mathematics.\footnote{Sections 3.7 to 3.12 in \cite{Vopenka:Alfa}
sketch the transformation of mathematics into set mathematics and subsequently set mathematics in Cantor's set theory. 
This sketch does mention ``various set theories'' but offers no explanation whatsoever as to which theories might
be meant (p.~126); towards the close of the exposition Vopěnka simply arrives at the following statement:
``The causes of the contemporary crisis of infinitary mathematics needs to be sought not in its manifestations,
but  in the  very foundations of the set-theoretic conception of mathematics. 
But this foundation is Cantor's set theory.'' (p.~128)}

The aforementioned dichotomy appears partly as a result of the \emph{$\atm$  narrative},
where the theory aimed to (in fact, in the \emph{Nová infinitní matematika} incarnation, it still does) 
compete with the extant---and inadequate---foundations. Prior to the $\atm$, 
a similarly great narrative featured the theory of semisets: ``The theory of semisets 
is not an auxiliary technical means but it studies a notion which might prove to be one
of fundamental mathematical (and metamathematical) notions.'' \cite{Hajek:whySemisets}
The simplest version of the narrative renders any other players, that is, other ambitious alternative
foundations of mathematics or indeed reasons why such alternatives ought to be taken into account,
as cameos. The dichotomous depiction is further underlined by a rather unfortunate name choice in English:
\emph{the} alternative set theory (this has been discussed in \cite{Chudacek:nekrologVopenka}).
Comparing expectations with the perceived reality, the narrative peters out in regret
that the $\atm$ has not been discovered, or been understood, or been appreciated, by the world at large. 

The dichotomy and the epic aspects of the $\atm$ can be documented in the literature;
we refrain from listing the references here, mainly because they are mainly  remarks or sidenotes in texts
aimed at discussing other topics. Thus they contribute to a depiction that surpasses any particular
text and persists as a tacit assumption. As one example, the author can offer an excerpt from her
own text \cite{Hanikova:HajekBiography}: 
``While Vopěnka’s alternative set theory is still a popular concept among Czech logicians, from
a more global point of view it seems to have shared the fate of many hitherto proposed alternatives 
to the mainstream conception of mathematics: it was trampled underfoot the crowd
that pursued the classical direction.” A text that aims to provide a scientific biography of P.~Hájek,
and thus it concerns P.~Vopěnka's work only indirectly, tends to slip into a simplified view 
on the standing of Vopěnka's $\atm$ school; that view is accepted implicitly.
 
Our sketch of a dichotomous presentation of the $\atm$ is in itself a simplification.
To be sure, there is no obligation for a working mathematician, seeking to free her/himself from the 
snares of the dominant mathematical ontology, to present an exhaustive survey of everything that has
ever happened outside this ontology. It may well be that in  opting for his method of presentation that
juxtaposes the classical set theory with the $\atm$, Vopěnka takes his cue from an assumed familiarity, or 
even preferences, of his readers. 

On the other hand, any work that seeks to reflect on Vopěnka's creative effort ought to place it in appropriate context.
Our view of the $\atm$ will be far from the dichotomous one: for context, we will consider not just the classical set theory,
but a synchronous and diachronous view of relevant parts of mathematics, including those that might be termed nonclassical.
We will highlight that the $\atm$ is not the singularity it is sometimes taken for, but rather
a part of well discernible currents in mathematics, namely the investigation of nonstandard universes (in Skolem's sense), strict
finitism in the sense of Parikh's paper \cite{Parikh:Feasibility}, and the effort to dislodge the artifacts of axiomatic systems
of formal theories\footnote{``It is commonly noted that set theory produces far more superstructure than is needed to
support classical mathematics.'' \cite{Holmes:SEPsettheory}.}, 
which manifests itself primarily in reverse mathematics. 
Our intent thus complements the historical facts already commonly listed by Vopěnka with respect to the $\atm$:
an analysis of Bolzano's \emph{Paradoxes of the Infinite}, Leibniz's infinitesimal calculus, 
occasionally selected citations of the French mathematical school of the turn of the 19th and 20th century, and
Robinson's nonstandard analysis.

An affinity between the $\atm$ and Robinson's nonstandard analysis is often considered. 
Vopěnka personally makes the following remark \cite[p.3]{Vopenka:Teubner}: ``From the formal and technical point of 
view, Alternative Set Theory is rather near to Nonstandard Analysis and can be considered, 
from this point of view, for a particular case of Nonstandard Analysis.'' 
We will put the main objection raised by Vopěnka against nonstandard analysis into
a context of his analogous objections against other areas of nonclassical mathematics, such
as the theory of fuzzy sets or intuitionistic logic: namely their \emph{subordination} to Cantor's set theory,
to which Vopěnka juxtaposes  a departure (ontological, and thereafter axiomatic) from this theory.
We will also recall that Vopěnka's work with nonstandard models of set theory starts approximately 
in 1960; before he could have had the opportunity to become acquainted with Robinson's publications on 
nonstandard analysis \cite{Robinson:NSA-1961,Robinson:NSAbook}.\footnote{This has been noted 
in publications concerning Vopěnka's work: \cite{Sochor:Vopenka, Svejdar:infinitenumbers}.}
 
To the extent in which works on the $\atm$ were published in English, 
they have been interpreted, placed in context, evaluated and applied by English-speaking researchers.
For an example, see Fletcher's chapter ``Infinity'' \cite{Fletcher:Infinity},
which takes the $\atm$ of Vopěnka and Sochor for an example of sophisticated theory of feasibility, or
Dean's paper ``Strict finitism, feasibility, and the sorites'' \cite{Dean:strictFinitismFeasibilitySorites}, where
the $\atm$ is cited as a presumably singular theory that uses nonstandardness to model vagueness.  
Recent mathematical work on the $\atm$ seem to appear rarely: let us mention Jeřábek's work ``Provability Logic of the
Alternative Set Theory'' \cite{Jerabek:thesisAST}.
Thanks to M.~R.~Holmes, an encyclopedic entry \cite{Holmes:SEPsettheory} on alternative axiomatic set theories
is available, and on the same theme, there is the paper with T.~Forster and T.~Libert \cite{HolmesForsterLibert:AST}.
It is therefore not correct to assert that Vopěnka's $\atm$ is forgotten. 
Along with his work in classical set theory, 
mainly boolean models and large cardinals, it has found a place within a collective awareness.

\section{On the mathematical design of the $\atm$\protect\footnote{
The sources for the present section are both of the already
mentioned monographs by Vopěnka on the $\atm$:  \cite{Vopenka:Teubner} and \cite{Vopenka:Alfa}.}}
\label{sct:matatm}
 
The $\atm$ provides self-contained foundations of mathematics without relying on the support of 
Zermelo--Fraenkel's set theory with the axiom of choice ($\zfc$) or any other 
already developed theory. The account of the mathematical world, supported by  suitably chosen axioms, 
is thus possible without reference to the ontology provided by classical set theory,
and by extension, to models of the $\atm$ in $\zfc$.\footnote{Not merely possible: 
the major part of \cite{Vopenka:Alfa} indeed proceeds in such manner.} 
One might say that an axiomatic nonclassical mathematics is developed over classical logic, 
remarking only that the
$\atm$  is sometimes presented as an informal theory, determined by its principles rather than
a collection of axioms. This maxim notwithstanding, several variants of an axiomatic
theory are available, first  outlined in Vopěnka's monograph \cite{Vopenka:Teubner}, 
while Sochor's follow-up works \cite{Sochor:ASTmeta-I,Sochor:ASTmeta-II,Sochor:ASTmeta-III} 
offer a precise analysis of the formulations of each of the axioms, their variants, and their interaction.
 
The monograph \cite{Vopenka:Teubner} presents the proof of each of the axioms of $\zffin$---a variant
of Zermelo-Fraenkel theory in which the axiom of infinity is replaced with its negation---in the $\atm$,
hence the full apparatus of the theory of hereditarily finite sets is available\footnote{
Sochor's paper \cite{Sochor:ASTmeta-III}, $\S$ 9 on models of the $\atm$ in $\zf$
shows that the $\atm$ is conservative over $\zffin$. 
This means that any statement in the language of sets is provable in $\zffin$ if and only if it is provable
in the $\atm$. An earlier 
$\S$ 6 of \cite{Sochor:ASTmeta-III} shows that the $\atm$ is equiconsistent with third-order arithmetic.} 
 for example, one can introduce the class $\N$ of natural numbers.
 
The $\atm$ uses classes as a standalone sort of objects.
It is only to be expected that classes are going to be an essential language tool for formulating
some interesting theses that are quintessential to the $\atm$.
The $\atm$ does not limit itself to set-definable classes.
\emph{Semisets} are subclasses of sets; proper semisets (those that are not sets themselves)
are not set definable.

An example of a class that is not set definable is the class $\Fin$ of sets
that are finite in the sense of $\atm$: that is, those sets \emph{every subclass of which is a set}. 
This definition gives rise to the notion of \emph{natural infinity}, the very \emph{raison d'{\^ e}tre}
of the $\atm$. All proper classes are (in Vopěnka's sense, i.e., naturally) infinite: 
``on encountering such a proper class, one always encounters the phenomenon of imprecision\footnote{Translator's note: 
the original Slovak term is ``neostrosť''. Rendering it here as ``imprecision'' is guided by 
\cite[p.~6]{Vopenka:Teubner}: ``Repeated generalization and higher precision have resulted 
in the determination of the form of the properties of objects so that the collection of all 
objects having such a property is understood as an actually existing set.''}
on the number of its elements [\dots] Natural infinity thus manifests itself, first and foremost, on proper classes.''
Notice that the notion of infinity is, in its definition, already linked to the notion of imprecision. 
A set that subsumes a proper class  is infinite in the sense of the above definition: 
``one encounters the phenomenon of imprecision even on the number of elements of some sets,
namely in case that, on encountering such a set, 
an imprecise part separates itself from it.'' (both quotations
\cite[pp.137--38]{Vopenka:Alfa}).
The essence of Vopěnka's infinity is captured quite well by A.~Sochor \cite{Sochor:BasesAST}: 
``The classical set theory puts infinity only `beyond' finite sets
while the alternative set theory puts it also among formally finite sets.''\footnote{ ``III. Every countable function 
is a subclass of a function which is a set. This thesis can be held for
an expression of the human aspiration of going beyond horizons---to transcend them mentally. 
It appears to be the most important and powerful thesis of the alternative set theory.
General mathematics works with infinity, therefore, admitting only formally finite sets, we
must be able to reach infinity through finite sets. Our third thesis expresses exactly a formalization 
of this approach (which is, in fact, the only possible in our case). Moreover, it makes [it]
possible to pass from investigation of infinite structures to investigations of formally finite
sets with convenient structures. The classical set theory puts infinity only `beyond' finite sets
while the alternative set theory puts it also among formally finite sets, since there are sets
containing infinite subclasses.'' \cite{Sochor:BasesAST}}

The intersection of the classes $\N$ and $\Fin$ is the class $\FN$ 
of all finite natural numbers; a flagship proper semiset. 
The prolongation axiom says that a (class) function on $\FN$ 
is a subclass of some set function. This implies the existence of
proper semisets: while the theory of semisets admits the existence of proper semisets, 
the $\atm$ postulates it. Within the $\atm$, the class $\FN$ is traditionally understood as
a semiset representing a \emph{path to the horizon}. Both $\N$ and $\FN$ yield interpretations of 
Peano arithmetic and hence they offer a basis for constructing nonstandard extensions of 
rational numbers. Infinite natural number are sets, and thus they obey, as 
sets, the laws of $\zffin$: no set bijection exists between two distinct natural numbers,
but on the other hand the class superstructure allows for a variant of 
 Hilbert's hotel, namely the construction of a bijective class map
of $x$ onto $x \cup\{ y\}$ for any infinite set $x$ and any set $y$.

An informal rendering of the train of thought expressed by the above constructions 
can be found in the review of \emph{The Theory of Semisets}  by A.~Lévy \cite{Levy:semisetReview}:
``For a particular model of set theory, there may be other collections of sets of the model
[than those that are deemed too big; author's note]
that cannot be admitted as sets of the model because they could cause the model to disobey
the axioms of set theory. A simple example is a collection of natural numbers without a least
member in a non-standard model of set theory. Such collections `know too much' even though
they are not necessarily very big. Conventional set theory does not admit the existence of such
collections either as sets or as classes. There is no way to speak about them in the theory---they
are `nonpersons'. The main idea of the book is to admit such collections as additional classes.''

The $\atm$ only uses two infinite cardinalities of classes, namely countable classes (which are proper semisets) 
and uncountable ones. 
The existence of a class bijection between any two uncountable classes
is postulated as an axiom and presented as the simplest position one can take 
as to cardinalities of classes, and as the maximal one regarding the existence of 
(bijective, class) maps \cite[p.~52]{Vopenka:Teubner}.

\section{The genesis of the $\atm$}
\label{sct:geneze}

In 1962 Vopěnka publishes a collection of works on the construction of a nonstandard
model of set theory\footnote{Starting with the paper \cite{Vopenka:OdinMetod}.},
where he extends Skolem's construction of a nonstandard model of Peano arithmetic \cite{Skolem:NonstandardPA}
which L.~S.~Rieger first introduced him to.\footnote{See \emph{Prolegomena} to \cite{Vopenka:NIM}, p.~26.} 
In his \cite{Vopenka:PragueSetTheorySeminar}, Vopěnka writes explicitly:
``At that time, I also began to frequent Rieger's seminar in mathematical logic, where I obtained
the knowledge of Skolem's nonstandard model of arithmetic of natural numbers. Thus armed, I created a nonstandard
model of G\"odel--Bernays set theory (that is, a model with nonstandard natural numbers) in 1961, using
the ultraproduct method.'' His work does not explicitly refer to Skolem's work, 
nor does it mention, e.g.,  {\L}o{\'s}'s ultraproduct theorem from 1955 \cite{Los:QuelquesRemarques};
in this case, Vopěnka confirmed many years later that his works are independent \cite[p.~724]{Vopenka:vyprNovobarokniMat}:
``I was not familiar with this work [{\L}o{\'s}'s paper; author's note] and neither were Tarski's students in the United States,
who worked on similar topics as I did at quite the same time. (At least Dana Scott, in the work discussed in 
the following section, does not refer to it.)''
A.~Sochor, in his paper \cite{Sochor:Vopenka}, added an exclamation mark to the year when discussing these works 
of Vopěnka from the early 1960s, which can perhaps be understood as relating to Robinson's
nonstandard analysis \cite{Robinson:NSA-1961,Robinson:NSAbook}, of which the given work of Vopěnka is independent;
Skolem's 1934 work is the common origin.

In 1972 Vopěnka and Hájek published, simultaneously with
\emph{Academia} and with \emph{North-Holland Publishing Company}, the monograph 
\emph{The Theory of Semisets} \cite{Vopenka-Hajek:Semisets}. 
The theory of semisets admits the existence of proper semisets, 
i.e., proper subclasses of sets.
Working with nonstandard domains provided  one of the views of universes of sets that 
 led to the introduction of the notion of semiset and its corresponding axiomatics: 
``There are also other axioms on semisets which extend the theory of semisets conservatively
w.r.t. statements which concern sets alone. For example, both the axiom of standardness and
its negation extend the theory of semisets in this way. The axiom of standardness asserts that
every non-empty semiset of ordinals has a first element [\dots] The theory of semisets with the
negation of the axiom of standardness will not be studied in the present work, but it can be
used in the study of non-standard analysis for example.” \cite[p.~12]{Vopenka-Hajek:Semisets}  
Other such domains are, for example,
generic extensions or, on the other hand, the subuniverse of constructible sets.
The theory of semisets aimed for very  high goals in the foundations of mathematics:
``In presenting the theory of semisets the authors hope to make some contribution to the task
of breaking through the bars of the prison in which mathematicians find themselves. The prison 
is set theory and the authors believe that mathematicians will escape from it just as they
escaped from the prison of three-dimensional space.'' \cite[p.~12]{Vopenka-Hajek:Semisets}  

The already mentioned set theory seminar at the Faculty of Mathematics and Physics of Charles University in Prague,
sometimes referred to as (Vopěnka's) Prague School of Set Theory, was a successor to a seminar that was run in Prague
by L.~S.~Rieger; after his death in 1963, the seminar run by Vopěnka continued until 1968, with outstanding results
some of which have become canonical within classical set theory.
The book \emph{The theory of Semisets} presents them in a semiset paradigm.
Lévy's review \cite{Levy:semisetReview}, from which we have borrowed the guideline on how to understand semisets,
ultimately views this fusion of results of the Prague School of Set Theory with
the semiset framework of presentation as the main weakness of the enterprise.

Robinson's nonstandard analysis stands out among the sources that Vopěnka's
work on the $\atm$ relies on, not least because he both refers to it and rejects
the way in which it is presented, which one might perhaps call semantic. Robinson works with a nonstandard structure and
its standard (elementary) substructure; in his 1961 paper there is a nonstandard extension 
of the structure of the real numbers, the 1966 book extends the already introduced techniques
to other structures; cf.~\cite{Robinson:NSA-1961,Robinson:NSAbook}.
Let us complement the quotation from our introductory section with another one, providing more information,
from \cite[p.~128]{Vopenka:Alfa}: ``The author of the alternative set theory admits 
that he, too, was encouraged  in creating this theory  by Robinson's nonstandard analysis,
through which the skeleton of this theory can be modeled in Cantor's set theory. 
At the same time, however, he has gained the experience that whenever 
he yielded to the temptation  to imitate the methods usual in nonstandard analysis, 
he was invariably drawn into the beaten path of mathematics based on Cantor's
set theory, which led him in undesirable directions, against the spirit of the alternative set theory.''
In contrast, Vopěnka's presentation is axiomatic; 
the following two excerpts offer insights into the genesis of the $\atm$ in the context of
Robinson's nonstandard analysis. The first comes from Hájek's article
``Why Semisets'' \cite{Hajek:whySemisets} from 1973: 
``But the use of semisets described by semiset theoretical axiom systems not consistently extendible 
to Set theory consists of a description and investigation of things (objects, situations)
set-theoretically not available. (i) [\dots] (ii) Various forms of axiomatic non-standard analysis.
(iii) [\dots]''; and in a footnote: ``Vopěnka; not published, but Vopěnka presented his
attempts in several lectures on various occasions.''
Thus, according to Hájek, Vopěnka attempts to axiomatize nonstandard analysis. 
In a similar vein, but  four decades later, P.~Pudlák \cite[p.~237]{Pudlak:LogicalFoundations}  remarks: 
``[\dots] Vopěnka proposed to use nonstandard analysis as the foundations of mathematics.
Instead of studying nonstandard models in Zermelo--Fraenkel theory, he suggested stating
the principles of nonstandard analysis as axioms and working in such an axiomatic system.''
V.~Švejdar, in his paper \cite{Svejdar:VopenkaHajek} gives an even more general statement:
``Vopěnka and Hájek believed that we (should) have the freedom to work with 
abstract axiomatic theories.''
 
Two more nonstandard set theories  are analyzed in the monograph
\emph{Nonstandard Analysis, Axiomatically} \cite{KanoveiReeken:NAA}. These are
Internal Set Theory ($\ist$)  as presented by Nelson \cite{Nelson:InternalSTpaper}
and Nonstandard Set Theory as presented by Hrbáček \cite{Hrbacek:NST},
which provided a basis for the cited monograph.
Both these theories start from $\zfc$, whereby they differ substantially from the $\atm$
(since the latter extends $\zffin$); the term ``nonstandard analysis'' is loosely applied
to indicate work with explicitly nonstandard universes. 
The authors remember Vopěnka's $\atm$ in the introduction,
where the sequence of the respective inceptions of the three theories is discussed \cite[p.~VII]{KanoveiReeken:NAA}:
``These three initial attempts to fully axiomatize nonstandard mathematics can hardly be 
linearly ordered in any reasonable sense, with any sort of preference assigned in some sound
manner. It is fair to assert, on the basis of available records, that all three were undertaken
independently of each other and led to results of comparable quality (although not of comparable 
impact on the practice of nonstandard mathematics, where IST has preference), in addition 
all three were based upon earlier development in foundations of nonstandard analysis.''
Yet another example of a nonstandard set theory is presented in Fletcher's paper \cite{Fletcher:NST},
which in turn reflects the earlier papers of Hrbáček and Nelson.

Another influence that can be discerned in Vopěnka's work on the $\atm$, 
even if it is not mentioned explicitly, is Parikh's paper on feasible arithmetic
\cite{Parikh:Feasibility}, namely his notion of \emph{almost consistent theories}. 
Speaking roughly, these are [consistent] extensions of Peano arithmetic with primitive recursive
functions which, in addition, postulate the existence of a cut (a unary feasibility predicate)
such that the value of a primitive recursive variable-free term $t$  does not
fall within the cut.
The theory is inconsistent; the paper shows that no contradiction can be obtained
from proofs of a bounded length. More precisely, 
a lower bound on the value of $t$ can be obtained as the value of a certain primitive recursive
function applied to any given bound on proof length,
and then proofs up to this length will not provide a contradiction when assuming
that $t$ is not feasible.\footnote{The Czech original  omits
by mistake the assumption that the theory without the feasibility predicate be consistent,
which is necessary for the result to hold as stated.} 
This material provided by Parikh was available to the group of Vopěnka's collaborators:
in Hájek's paper \cite{Hajek:whySemisets} from 1973 the collection of feasible numbers
is rendered as a semiset. Vopěnka's work occasionally mentions the notion of 
a witnessed universe: one in which a proper semiset can be found within a concrete set\footnote{see \cite[p.~37]{Vopenka:Teubner}}.
One may wonder what is meant by a \emph{concrete} set.
But the same passage of text asks the reader to accept that large standard numbers,
represented as values of primitive recursive functions (e.g., $67^{293^{159}}$)
might, as sets, subsume a proper semiset, and hence might be infinite in the $\atm$ sense.
As Vopěnka admits, such ideas are not elaborated mathematically;
axiomatic $\atm$ can be understood as a ``limit case'' of such considerations\footnote{
``The axiom of existence of proper semisets does not imply that there are proper semisets
included in any specific set. If we guarantee the existence of a proper semiset included in
a certain concrete set then we say we are studying a witnessed universe. When we restrict
the family of properties admitted in the axiom of class existence in such a way that proper
subsemisets of all concrete sets are eliminated we say that we are studying a limit universe.''
\cite[p.~37]{Vopenka:Teubner}}.

Vopěnka's work in the $\atm$ is sometimes mentioned as modelling feasibility by means of nonstandardness (see, e.g., \cite{Fletcher:Infinity}),
with the tendency to identify the meaning of ``standardness'' with that of ``feasibility''
and to view the $\atm$ as a (consistent) mathematical model of the feasibility notion,
without touching upon the problem of unfeasibility of large standard natural numbers.
Other authors view Vopěnka's work in the $\atm$ as an attempt to model vagueness (see, e.g., 
\cite{Bellotti:reductionInfinite,Dean:strictFinitismFeasibilitySorites}).

\section{A context for the $\atm$}

Our analysis of the mathematical origins of the $\atm$ aims to bring a few conclusions.
Firstly, the theories, notions, and results discussed above are, in themselves, contributions
to the map of the world of mathematics, which we prefer over the dichotomous view.
Mathematically, the $\atm$ extends Vopěnka's active, long-term experience with the results
in classical mathematics, or possibly with selected areas of nonclassical mathematics, if we
choose to view, e.g., nonstandard universes or Parikh's work on almost consistent theories as nonclassical.\footnote{
Parikh's paper \cite{Parikh:Feasibility} provided an impulse for the study of bounded arithmetic, cf.~\cite{Buss:TwoPapersParikh}.} 
In particular: the $\atm$ uses the notion of semiset, introduced in the monograph \cite{Vopenka-Hajek:Semisets},
it processes Vopěnka's experience with work in nonstandard models of set theory, and it brings his own axiomatic conception
of nonstandard analysis.

In the period that was key to the inception of the $\atm$, i.e., the turn of the 1960s and the 1970s,
Vopěnka was prevented from publishing his work on political grounds, according to his contemporaries' recollections,
cf.~\cite{Chudacek:nekrologVopenka}. There is a marked absence of primary sources on the $\atm$, which one
might expect under normal circumstances and which would have documented its genesis.
From this point of view, Hájek's paper \cite{Hajek:whySemisets} from 1973 seems rather valuable,
given that it places Vopěnka's attempts to axiomatize Robinson's nonstandard analysis into this period.

Secondly, without a further discussion one cannot accept the thesis that the $\atm$ was the outcome
of Vopěnka's effort at a phenomenological catharsis of mathematics. Various remarks
on philosophical motivations of the $\atm$ are impossible to miss when studying the primary and
secondary literature. Let us quote, for example, P.~Zlatoš in \emph{Filosofický časopis}:
``I could not bring myself to see how, on the basis of such incomprehensible, intimidating 
treatises [the works of Husserl and Heidegger; author's note], Vopěnka could have arrived at his
beautiful, crystalline reflections on the horizon and its role in our understanding of the world.'' 
\cite[p.~520]{Zlatos:ATMUdelGenia} or ``What ought to be considered the permanent 
contribution of the $\atm$ and of the philosophical considerations
it originated from?'' \cite[p.~531]{Zlatos:ATMUdelGenia} 
Another example can be found in the text of A.~Vencovská from the same issue:
``A key criterion for the correctness of any reasoning on infinitely small quantities in nonstandard analysis is
whether or not they behave as the entities in the ultraproduct that represent them. 
This is the reason why the Cantorian intuition on infinity, as it is captured in classical set theory,
is still a key one in nonstandard analysis. 
On the other hand, Vopěnka's intuition originated from his philosophical considerations.
Even as he was aware of the influence of Robinson's work on himself, he obtained the principles
that he adopted for his new infinitary mathematics from his philosophy.'' \cite[p.~585]{Vencovska:MnohoPovyku}
But it is presumably Vopěnka himself,  albeit  in retrospect,
who is the staunchest proponent of the thesis. For example, his essay \cite{Vopenka:ASTallabout}\footnote{
The English translation was provided by A.~Vencovská. Later published as \cite{Vopenka:AstPhil1991};
see also \cite{Vopenka:ocojdevatm}. } 
offers the following: ``In my book \emph{Mathematics of the Alternative Set Theory} I overestimated
the readers. I believed that the philosophy of the Alternative Set Theory would be immediately obvious to all
and I therefore concentrated more on demonstrating some techniques. However, the traditional mathematical thinking is so 
established---even in these times when physics is reevaluating its traditional thinking---that most mathematicians
are not able to carry out a phenomenological critique of classical ideas, let alone realize its consequences.
Therefore in this paper I shall attempt to sketch, at least briefly, the philosophical foundations
of the Alternative Set Theory.''

From a phenomenological critique of the foundations of mathematics, there is no straight
  route to Vopěnka's design of the $\atm$ and its mathematization of key notions, such as that of infinity.
Vopěnka (with collaborators) discovered, developed, and established the affinity between
the Central-European phenomenological tradition and the mathematical methods of capturing nonstandard universes axiomatically.
These two influences blended happily in that particular time and place. Our paper does not purport to 
belittle the importance of this discovery, nor of Vopěnka's crusade for it, 
face to face with the initial void that is a typical attribute of free and creative research work.
Nevertheless, one needs to bear in mind that the nuptials of the two influences were very much an ongoing process,
and that Vopěnka entered into it equipped with approximately a decade of active, excellent research work in
set theory, and in particular, active work in the theory of nonstandard universes.

Philosophical context for the $\atm$ can be viewed as a \emph{validation}
of the available mathematical developments. Let us take the following lines from \emph{Mathematics 
in the Alternative Set Theory} as an example:
``We shall not always use all properties of indiscernibility relations and therefore we could
study more general topologies as in Cantor's set theory. The author has developed general
topology in the alternative set theory; interesting results in such topology were obtained by
J. Chudáček. But we shall not develop general topology here before deeper motivations for its
study in alternative set theory are exhibited.'' \cite[p.~87]{Vopenka:Teubner} 
Another example can be found earlier in the monograph:
``Our theory offers various possibilities of theories of infinity. For example, we could imitate
Cantor's whole theory. [\dots] We can assume axioms for any other theory of infinity, provided
it does not contradict the other axioms. Cantor’s theory is just one possibility. At present, no
reasons for the acceptance of a nontrivial theory of infinity are known. All such theories must
be speculative in character. Consequently, their results mentioning infinite cardinalities will be
vacuous if their speculative background is rejected. To prevent this, we decide to accept a 
trivial theory of infinite cardinalities.'' \cite[pp.~51--52]{Vopenka:Teubner} This text introduces the axiom of two cardinalities.
These samples suggest that a Vopěnka would not consider a mathematical treatment of a topic satisfactory
unless it fit in with the conception that apparently evolved in parallel to the mathematical considerations.
In the introduction to the monograph, Vopěnka writes: ``My colleague J.~Polívka helped me to improve various 
motivations.'' \cite[p.~3]{Vopenka:Teubner}
J.~Polívka was an external member of the Department of Mathematical Logic at the Faculty of Mathematics
and Physics of Charles University from 1969 till 1980 and a student of J.~Patočka since 1968.
As an example of a validation of a mathematical development with a philosophical thesis one can take
supporting the axiom of prolongation with the principle of surpassing the horizon.  

A reminder that possible phenomenological motivations of the $\atm$ may have been
more or less \emph{ex post} is present also in the reminiscences given in Vopěnka's obituary by Chudáček \cite{Chudacek:nekrologVopenka},
who was Vopěnka's PhD student between 1972 and 1976:
``Vopěnka provided a prolongation beyond the horizon in mathematics, similar to Husserl, although,
as far as I know, he discovered it independently of Husserl. He learned about Husserl and his work 
only later in the 1970s from the philosopher Jiří Polívka.''

Let us recall two important research areas that were close in spirit to Vopěnka's $\atm$ even
though there was no interaction. 
There are important parts of Vopěnka's essay ``What is the Alternative Set Theory all about'' \cite{Vopenka:ASTallabout} 
that concern  the contribution of Vopěnka's conception of infinity to the study of vagueness.
``Natural infinity is an abstraction from path leading to the boundaries of indefinite
and blurred phenomena. Since we understand various phenomena exactly through
grasping their boundaries (compare lat.~`definitio'), the study of natural infinity is
also a possible foundation for the science of indefiniteness.
[\dots]
That basic form of natural infinity which led, after completing the sharpening
process, to the basic classical infinity, countability, has been abstracted from the path
leading on to the horizon.''
Vopěnka holds that classical mathematics, on the other hand, forwent the challenge to model the phenomenon of indefiniteness.
Approximately at the same time when Vopěnka  writes his \emph{Mathematics in the Alternative Set Theory},
M.~Dummett publishes the essay ``Wang's paradox'' \cite{Dummett:WangParadox} that,
according to Dean's paper  \cite{Dean:strictFinitismFeasibilitySorites}, is
considered the origin of  contemporary philosophical study of vagueness (and it is also 
a detailed analysis---and a rejection---of ultrafinitism).

At the same time, H.~Friedman publishes first works on \emph{reverse mathematics}.\footnote{See,
e.g., \cite{Friedman:SysSecondOrderArithmetic}.}
The main idea is to seek suitable axioms for proving individual mathematical statements, 
where as a verification that the axioms were correctly identified, one would prove the axioms from the statements.
This project was mainly carried out in fragments of second-order arithmetic.
The theory is strong enough to formalize a major part of everyday mathematics.
The amount of results that reverse mathematics can boast may give the impression that
it has almost been ``done''. It is not in a close relation to the $\atm$, but  it has
established itself as a serious attempt to get rid of the artifacts of axiomatic theories: and it does
carry out this intent, cherished by Vopěnka, very successfully.  
Nevertheless, the framework of reverse mathematics has not been quite accepted as a serious
candidate for new foundations of mathematics; admittedly it could play this role, but it does not play it.

\section{The $\atm$ and nonclassical mathematics}

We only have fragments of evidence as to Vopěnka's views on the work of his companions 
within the wide plains of nonclassical mathematics. What stands out is Robinson's 
nonstandard analysis, often mention with a shade of reservation.
The book \emph{Mathematics in the Alternative Set Theory} is one among many
works that only touch upon the publications of Yessenin-Volpin.\footnote{
Parikh's paper ``Existence and Feasibility in Arithmetic''  \cite{Parikh:Feasibility} belongs in this group too. On the other hand,
Dummett's ``Wang's Paradox'' \cite{Dummett:WangParadox} 
does go into quite some depth analyzing ultrafinitism as represented by Yessenin-Volpin.}
The monograph \emph{Vyprávění o kráse novobarokní matematiky} [The Tale of the Beauty of Neobaroque Mathematics]
\cite{Vopenka:vyprNovobarokniMat} brings its own survey of the history of set theory,
including some less than orthodox branches, such as the axiom of determinacy.

A departure from classical mathematics can be induced already at the level of logic
in which a mathematical theory  is developed.
Separating a ``logic'' from an ``axiomatic theory'' may appear way too formalistic,
however one can indeed trace such a view in Vopěnka's work. 
The $\atm$ lives in classical logic; Vopěnka may have intended an unspecified departure from 
classical logic, but it was not enough of a priority to him to carry it out explicitly.
An early mention of it can be found in \cite[p.~11]{Vopenka:Teubner}:
``The careful reader will realize that our theory is not even fully based on the classical concept
of finite sets. This is only implicit in the present book since systematic development in this direction 
would make the book considerably longer and would stress topics not suitable for first
acquaintance with our theory. We have the same reason for development of our theory inside
classical logic.''
Another faint hint can be found in \emph{Nová infinitní matematika} \cite{Vopenka:NIM}:
``When studying the semiset $\FN$, its properties and their mutual relationships therefore, we  can
use the predicate calculus. It is however necessary to use quantifiers with caution, since
towards the horizon that bounds the sizes of finite natural numbers, the semiset $\FN$ is not sharply delimited.''
[\emph{Prolegomena}, p.~41]

Another point of interest lies in Vopěnka's position on intuitionism and intuitionistic logic.\footnote{See
also remarks on this in the last part of Švejdar's paper \cite{Svejdar:VopenkaHajek}.}
Vopěnka viewed it as a restriction of classical logic as to the available techniques.
In \cite[p.~103]{Vopenka:Alfa} he gives the following assessment:
``A radical way out of the crisis [occasioned by laying out  contradictions and
paradoxes in naive set theory; author's note] was offered by Brouwer's intuitionism,
which completely abandoned the actualization of infinite domains of objects as unjustified,
whereby it parted ways entirely not only with Cantor's set theory, but also with
the general part of Bolzano's account of infinity. [\dots]
Unfortunately, such a revision of logic is not compatible with the main intent of set theory,
which aims precisely at the actualization of infinite domains of objects.''
Within the context of classical set theory, intuitionism is presented as a dead end.
Many years later, in \cite[s.~698--700]{Vopenka:vyprNovobarokniMat}, Vopěnka writes:
``As soon as the intuitionists clearly established the laws of logic that they were using,
it was only to be expected that these would be noticed even by those mathematicians who
tried to interpret mathematics as a subarea of set theory.
Indeed, it was not long after the publication of the Heyting work we mentioned that a surprising
colonization of intuitionistic logic by the set-theoretic empire took place. [\dots]
Stone's truth values of propositions in intuitionistic logic present a grand victory
of the set-theoretic empire over intuitionism.''
These comments are analogous in nature to those raised with respect to nonstandard analysis:
there is a semantics here that is presented within classical set theory.
But intuitionism, just like the $\atm$, originates from principles, whereupon an axiomatization follows.
From the perspective of this parallel, it is not clear how a topological semantics of intuitionistic 
logic within classical set theory could constitute a serious objection against intuitionism.

Given the essential connection between indefiniteness and the $\atm$, let us 
briefly discuss fuzzy logic and fuzzy sets that provide a different modelling of indefiniteness.
It remains unclear whether Vopěnka was familiar with fuzzy logic as a formal nonclassical system 
(despite the fact that, incidentally, just at the same time the key publication on the $\atm$ was being written,
J.~Pavelka, a student of A.~Pultr, was writing his dissertation on fuzzy logic at the Faculty of Mathematics and Physics;
this work \cite{Pavelka:Fuzzy} has since become canonical).
The theory of fuzzy sets, in Vopěnka's view, fails to escape the clutches of the world of classical mathematics,
because it uses real-valued sets: ``As a matter of fact, classical mathematics only hamstrings
the progress of the theory of fuzzy sets, as the former leads the latter down a path where 
it is hard to understand why the latter might venture of its own will.'' \cite[p.~94]{Vopenka:Alfa} 
This objection, however, presumably conflates the theory with its model in the classical universe of sets.
Granted, Zadeh's pioneering paper \cite{Zadeh:1965} made popular the notion of fuzzy set taken as a function
that assigns to each object of a given universe a real number in the unit interval.
The semantics of reals also determines the set-theoretic operations. But even Zadeh's approach
admits general functions, alongside real-valued ones. This was continued in J.~A.~Goguen's essay
``The logic of inexact concepts'' \cite{Goguen:LogicInexactConcepts}, which in turn inspired Pavelka's work
on axiomatization of the required formalism. Many works on axiomatic set theory in fuzzy logic 
make hardly any use of classical logic and classical set theory;
as a representative, let us name ``Bemerkungen zum Komprehensionsaxiom“ \cite{Skolem:Komprehensionsaxiom} from 1957, 
authored once again by Th.~Skolem. 

These detours are intended further to illustrate how crucial  it was to Vopěnka 
that an \emph{axiomatic} theory be independent of any of its \emph{classical} interpretations.
Our mathematical reasoning needs to be supported by principles rather than by
an interpretation in classical set theory that may have been constructed to ``illustrate'' them;
going in the reverse order  would only render the new theories subordinate to the classical theory. 
Bellotti offers the following thesis: ``[\dots] the most natural interpretations of $\atm$ 
are in terms of non-standard models of arithmetic, so that one is strongly tempted to 
reduce this approach conceptually to some non-standard background, 
although Vopěnka would reply that this would reverse the conceptual priority, 
since non-standard analysis is for him a (rather artificial) `Cantorian' way to treat `natural’
infinity.'' \cite{Bellotti:reductionInfinite}
Still, both intuitionistic logic and intuitionistic set theory on the one hand, and
fuzzy logic and fuzzy set theory on the other, did in fact attain this kind
of independence from classical mathematics.

\section{Conclusion}

The $\atm$ has its place in the mosaic of the mathematical world: we
provide an alternative to the assumption that this world is too homogeneous and definite 
whereby the $\atm$ does not fit in,  keeps its distance and merely looks on.
We have described the mathematical substrate of the 20th century that the $\atm$ builds on
and have pointed out some parallel developments that share some of their starting points with
the $\atm$ and underline the plurality of the foundations of mathematics.
We have also recalled, in several instances,  Vopěnka's preference for axiomatic method of work.

The present text is about the mathematical neighbourhood of the $\atm$,
and the present author's maxim has been to adhere to sources available in writing.
Which is the inner structure of its philosophical neighbourhood? 
What is the genesis of the philosophical frame of the $\atm$ and how is it anchored in time?
Which terms would aptly describe the viable topics of this genesis,
and for each topic, which would be a suitable comparative commentary or
corresponding passages of other researchers' works?

The text ``What is the alternative set theory all about'' \cite{Vopenka:ASTallabout}
from 1989 is presumably among the most fruitful resources for a reflection of this kind.\footnote{Another
resource is the paper \cite{Trlifajova:NekonecnoVopenky} and the monograph \cite{Vopenka:Meditace} it often refers to.}
It indicates the philosophical foundations of the $\atm$ as its agenda and
offers, among other things, a detailed exposition of Vopěnka's natural infinity 
in contrast to the classical one.
But at the same time, it points out problems awaiting the reflecting researcher (for example,
it gives no references). The text also appears long after  the key period for the $\atm$,
which was the early 1970s. In the geology of the $\atm$ therefore, it acts as 
a smooth river rock, since to the inexperienced eye,  there is no telling 
as to which cliff it was once torn out of by the elements.
Yet there is still hope that experienced geologists could provide interesting accounts
of the journeys undertaken by this rock and perhaps other ones.

\bigskip
\noindent{\bf Acknowledgements}
The preparation of this paper was supported by the long-term strategic development financing
of the Institute of Computer Science (RVO:67985807).  Libor Běhounek, Jirka Hanika and Štěpán Holub 
provided valuable comments on some earlier versions of the Czech text or its parts. The ICS Library, particularly
Ludmila Nývltová, offered help in obtaining  older references. Thanks are due to three anonymous reviewers
appointed by \emph{Filosofický časopis} for their comments that improved the Czech text. 

This is the author's translation of her paper published originally in Czech: 
``Vopěnkova alternativní teorie množin v matematickém kánonu 20.~století'',
\emph{Filosofický časopis} 70(3), 485--504, 2022. {\tt https://doi.org/10.46854/fc.2022.3r.485}. 
The excerpts from references available only in Czech/Slovak were also translated by the author.
Drafts for the translation of several paragraphs in section 2 and 3 of the translation were first obtained using the
{\tt www.DeepL.com/Translator} (free version).

\bibliographystyle{plain} 

\renewcommand\refname{References}

\end{document}